\documentclass[10pt, a4paper, reqno]{amsart}
\usepackage{style}

\begin{document}

% Abstract
\begin{abstract}
We show that a formal Deligne--Mumford stack is formal-locally represented by a formal scheme.
This is an analogue of Frobenius theorem for smooth foliations in any characteristic and without smoothness hypotheses on the ambient space.
\end{abstract}

\maketitle

\tableofcontents
\section*{Introduction}
\label{sec:introduction}

In the past century, foliations have been intensively studied in the context of differential geometry.
However, they have recently received a surge of interest from algebraic geometers as a tool to tackle long standing conjectures, such as the abundance conjecture and the Green--Griffiths conjecture.
Following work of McQuillan (see \cite{MR4507257}), it seems clear that the correct analogue of a foliation in algebraic geometry is a formal groupoid.
This is a particularly well-behaved notion for studying singular foliations or foliations over fields of positive characteristic.
The aim of this article is to prove the analogue of Frobenius theorem for formal Deligne--Mumford groupoids (see Definition \ref{def:formal_dm_stack}).
These groupoids should be thought as the analogue of smooth foliations.

\begin{thmx}
\label{thm:frobenius}
    Let $R \rightrightarrows X$ be a formal Deligne--Mumford groupoid locally of formal finite presentation over a locally Noetherian formal scheme $S$ and let $x \in X$ be a closed point.
    Let $\hat{X}$ be the infinitesimal neighbourhood of $x \in X$ and let $R \hat{|}_{\hat{X}}$ denote the infinitesimal restriction of $R$ to $\hat{X}$.
    Then there exists a unique formal scheme $W$ over $S$ and a unique $S$-morphism $q : \hat{X} \rightarrow W$ such that
    \begin{align}
    \label{eq:frobenius_graph}
        \hat{X} \times_{W} \hat{X} = R \hat{|}_{\hat{X}} 
    \end{align}
    as groupoids on $\hat{X}$.
    Furthermore $q$ is formally smooth, split and locally of formal finite presentation.
\end{thmx}

Theorem \ref{thm:frobenius} can be seen as a consequence of a formal-local structure theorem for \emph{infinitesimal} Deligne--Mumford groupoids.
An infinitesimal groupoid is a formal groupoid which only identifies close by points (see Definition \ref{def:infinitesimal_groupoid}).
The next theorem states that if the topological space of an infinitesimal Deligne--Mumford stack is a single point, then the stack is represented by a formal scheme.

\begin{thmx}
\label{thm:local_structure}
    Let $R \rightrightarrows X$ be an infinitesimal Deligne--Mumford groupoid locally of formal finite presentation over a locally Noetherian formal scheme $S$ and suppose that $|X|$, the topological space of $X$, consists of a single closed point $x$.
    Then there exists a unique formal scheme $W$ over $S$ and a unique $S$-morphism $q : X \rightarrow W$ such that
    \begin{align}
    \label{eq:frobenius_local_graph}
        X \times_W X = R
    \end{align}
    as groupoids on $X$.
    In addition, $q$ is the categorical quotient of $X$ by $R$ in the category of formal schemes over $S$ and it represents the stack $[X/R]$.
    Furthermore $q$ is formally smooth, split and locally of formal finite presentation.
\end{thmx}

The proof consists of constructing a transversal (or slice) through the closed point as in \cite[Lemma 3.3]{MR1432041}.
This gives rise to a minimal presentation and, using the fact that the groupoid is infinitesimal, we can conclude that the transversal represents the associated stack.
The construction of the transversal is fairly technical as there is the additional difficulty of dealing with formal schemes. \par

A very similar statement appears in \cite[Fact 2.2]{MR4507257} where the author employs the same strategy.
However, we construct a transversal for formal schemes, rather than usual schemes, and we use different and more rigorous methodology. \par

We put special care into not using results for algebraic stacks.
In fact, we only work with presentations of algebraic stacks, i.e. groupoids.
This choice is justified by the fact that developing foundations for formal stacks requires some work, e.g.the choice of a suitable topology for the site and the proof of associated descent results. \par

Before giving a brief outline of the paper, we discuss why it is reasonable for formal groupoids to be the analogue notion of foliations.
A characterisation of a smooth foliation on a smooth manifold $M$ is an atlas of compatible equivalence relations on sufficiently small open sets.
This is the statement of the classical Frobenius theorem.
One can imagine that glueing all such equivalence relations will give a subset $R$ of $M \times M$.
This will not necessarily be an equivalence relation as it will fail to identify points in different charts, however it will be an \emph{infinitesimal equivalence relation}. \par

Algebraically, one can cook up an infinitesimal equivalence relation from a foliation $\mathscr{F}$ on a smooth variety $X$ in characteristic zero as follows: we can consider the sub-algebra $\mathrm{Diff}_X (\mathscr{F})$ of the sheaf of differential operators $\mathrm{Diff}_X (\mathscr{O}_X)$ of $X$ generated by $\mathscr{F} \subseteq \mathscr{T}_X$.
This is reasonable as the algebra of differential operators of a smooth variety in characteristic zero is generated by $\mathscr{T}_X$.
We can then dualise the inclusion $\mathrm{Diff}_X (\mathscr{F}) \subseteq \mathrm{Diff}_X (\mathscr{O}_X)$ so to obtain a closed formal subscheme $R$ of $( X \times X )_{\widehat{\Delta}}$, the completion of the diagonal morphism $\Delta : X \rightarrow X \times X$.
The formal scheme $R$ is the associated infinitesimal equivalence relation.
There are several details to check and a rigorous correspondence will appear in future work. \par

We use the conventions and results developed in \cite{bongiorno1}.
In particular, a formal scheme is a locally topologically ringed space which is locally isomorphic to the formal spectrum of an adic ring with finitely generated ideal of definition.

\subsection*{Outline}
In \S \ref{sec:infinitesimal_groupoids}, we recall the notion of formal and infinitesimal groupoids from \cite{MR4507257}.
In \S \ref{sec:infinitesimal_equivalence_relations}, we study infinitesimal equivalence relations and we show that an infinitesimal Deligne--Mumford groupoid is necessarily an equivalence relation.
In \S \ref{sec:restrictions}, we study restrictions and infinitesimal restrictions of groupoids to infinitesimal neighbourhoods of closed points.
In \S \ref{sec:transversals}, we define what a transversal is and show it always exists under Noetherian assumptions.
Finally, in \S \ref{sec:frobenius_theorem} we prove the main theorems.

\subsection*{Acknowledgements}
I would like to thank my PhD advisor Paolo Cascini for valuable advice.
I thank Lambert A'Campo, Alessio Bottini, Riccardo Carini, Przemys\l{l}aw Grabowski and Maarten Mol for answering my questions especially in Complex and Differential Geometry, and in positive characteristic.
This research was conducted at the Max-Planck-Institut f\"{u}r Mathematik and the Universiteit van Amsterdam, which I thank for their hospitality and financial support.

\section{Infinitesimal groupoids}
\label{sec:infinitesimal_groupoids}

A \emph{formal groupoid} $\mathcal{X}$ over a formal scheme $S$ is a groupoid object in the category of formal schemes over $S$.
This is the datum of two formal schemes $X$ and $R$ over $S$ together with five $S$-morphisms
\begin{enumerate}
    \item $s : R \rightarrow X$ (source);
    \item $t : R \rightarrow X$ (target);
    \item $e : X \rightarrow R$ (unit);
    \item $i : R \rightarrow R$ (inverse);
    \item $c : R \times_{(s,t)} R \rightarrow R$ (composition);
\end{enumerate}
such that for all formal schemes $T$ over $S$,
\begin{align}
\label{eq:functor_groupoid_description}
    \mathcal{X}(T) = \left( \mathrm{Hom}_S \left(T, R \right) \xrightrightarrows[t(T)]{s(T)} \mathrm{Hom}_S \left(T, X \right) \right).
\end{align}
is a groupoid of sets.
The morphism $j = t \times s : R \rightarrow X \times_S X$ is the \emph{diagonal} morphism of the groupoid.
A formal groupoid will also be denoted by $R \rightrightarrows X$. \par

If $s$, or equivalently $t$, is a morphism of formal schemes with a given property $P$, then we say that the groupoid has property $P$.
This will be applied with properties including smoothness, formal smoothness and formal finite presentation. \par

Since formal schemes are colimits of schemes, it suffices to check whether $\mathcal{X}(T)$ is a groupoid only for affine schemes.
In particular, this shows that usual groupoids of schemes are formal groupoids. \par

A formal scheme $X$ over $S$ can be viewed as a formal groupoid by letting $R = X$ and all structural morphisms to be the identity. \par

Let $\mathcal{X} = R \rightrightarrows X$ and $\mathcal{Y} = Q \rightrightarrows Y$ be formal groupoids over a formal scheme $S$.
A morphism of groupoid $f : \mathcal{X} \rightarrow \mathcal{Y}$ is a pair of $S$-morphisms $X \rightarrow Y$ and $R \rightarrow Q$ such that, for all formal schemes $T$ over $S$, the induced map $\mathcal{X}(T) \rightarrow \mathcal{Y}(T)$ is a functor of categories.
The category of formal schemes over $S$ embeds fully faithfully in the category of formal groupoids over $S$.
Note that there is a natural morphism of formal groupoids $X \rightarrow \mathcal{X}$.

\begin{lemma}
\label{lem:unit_immersion}
    Let $\mathcal{X} = R \rightrightarrows X$ be a formal groupoid over a formal scheme $S$.
    Then the unit morphism $e : X \rightarrow R$ is an immersion of formal schemes.
\end{lemma}

\begin{proof}
    There is a factorisation
    \begin{align}
    \label{eq:diagonal_factorisation}
        \Delta_{X/S} : X \xrightarrow{e} R \xrightarrow{j} X \times_S X,
    \end{align}
    where $\Delta_{X/S}$ is the diagonal morphism associated to the structural morphism $X \rightarrow S$.
    By \cite[Lemma 2.8]{bongiorno1}, $\Delta_{X/S}$ is an immersion of formal schemes.
    But then $e$, being the first factor of (\ref{eq:diagonal_factorisation}), is an immersion of formal schemes (\cite[Lemma 2.9]{bongiorno1}).
\end{proof}

\begin{definition}
\label{def:infinitesimal_groupoid}
    An \emph{infinitesimal groupoid} $\mathcal{X}$ over $S$ is a formal groupoid over $S$ such that $\mathcal{X}(T)$ is a trivial groupoid of sets (no non-identity morphisms) for all reduced schemes $T$.
\end{definition}

\begin{lemma}
\label{lem:infinitesimal_groupoid_characterisation}
Let $\mathcal{X}$ be an infinitesimal groupoid over a formal scheme $S$.
The following conditions are equivalent:
\begin{enumerate}
    \item $\mathcal{X}$ is an infinitesimal groupoid.
    \item The unit morphism is a thickening of formal schemes (see \cite[\S 3]{bongiorno1}).
    \item The source (or target) morphism induces a homeomorphism of topological spaces.
\end{enumerate}
\end{lemma}

\begin{proof}
    $(1) \rightarrow (2)$.
    We know already that $e : X \rightarrow R$ is an immersion of formal schemes (Lemma \ref{lem:unit_immersion}) and it suffices to show it induces a bijection of underlying sets.
    By assumption, whenever $T$ is reduced, the map of sets
    \begin{align}
    \label{eq:unit_bijection}
        e(T) : \mathrm{Hom}_S \left(T, X \right) \rightarrow \mathrm{Hom}_S \left(T, R \right)
    \end{align}
    is a bijection.
    Applying this when $T$ is the spectrum of a field yields that $e$ is a bijection. \par

    $(2) \rightarrow (1)$.
    Suppose $e$ is a thickening and let $T$ be a reduced scheme.
    Then any morphism $T \rightarrow R$ must factor through $e : X \rightarrow R$.
    This follows from observing that $e$ induces an isomorphism between the reduction of $X$ and the reduction of $R$.
    This shows that (\ref{eq:unit_bijection}) is surjective.
    Furthermore, since $e$ is an immersion, it is a monomorphism, hence $e(T)$ is injective.
    Therefore there are no non-identity morphisms in the induced groupoid of sets $\mathcal{X}(T)$. \par
    
    $(2) \leftrightarrow (3)$.
    We know that $e$ is a section of $s$.
    Therefore, $e$ induces a homeomorphism of topological spaces if and only if does $s$.
\end{proof}

Given a formal groupoid $\mathcal{X} = R \rightrightarrows X$, there is a natural infinitesimal groupoid $\hat{\mathcal{X}} = \hat{R} \rightrightarrows X$ associated to $\mathcal{X}$.
This is obtained by completing along the unit morphism $e$.

\begin{construction}
\label{cons:infinitesimal_from_formal}
    Let $\mathcal{X} = R \rightrightarrows X$ be a formal groupoid over a formal scheme $S$ and suppose $s : R \rightarrow X$ is locally of formal finite presentation.
    Then $e : X \rightarrow R$ is an immersion of formal schemes (Lemma \ref{lem:unit_immersion}) and is locally of formal finite presentation (\cite[Lemma 4.13]{bongiorno1}).
    Define $\hat{e} : X \rightarrow \hat{R}$ to be the infinitesimal neighbourhood of $e$ (see \cite[\S 5]{bongiorno1}) and let $\iota : \hat{R} \rightarrow R$ be the induced morphism satisfying $e = \iota \circ \hat{e}$.
    Firslty, we can easily define the source morphism $\hat{s} = s \circ \iota$ and the target morphism $\hat{t} = t \circ \iota$.
    We can define $\hat{i}$ by applying the universal property of infinitesimal neighbourhoods to the diagram
    \begin{equation}
    \label{diag:infinitesimal_inverse}
        \begin{tikzcd}
            X \arrow[d, "\mathds{1}"] \arrow[rr, "\hat{e}"] & & \hat{R} \arrow[d, "i \circ \iota"] \arrow[dl, dashed, "\hat{i}"] \\
            X \arrow[r, "\hat{e}"] & \hat{R} \arrow[r] & R .
        \end{tikzcd}
    \end{equation}
    Similarly, the morphism $\hat{c}$ exists by applying the universal property of infinitesimal neighbourhoods to the diagram
    \begin{equation}
    \label{diag:infinitesimal_composition}
        \begin{tikzcd}
            X \arrow[d, "\mathds{1}"] \arrow[rr, "\hat{e} \times \hat{e}"] & & \hat{R} \times_{(\hat{s}, \hat{t})} \hat{R} \arrow[d, "c \circ (\iota \times \iota)"] \arrow[dl, dashed, "\hat{c}"] \\
            X \arrow[r, "\hat{e}"] & \hat{R} \arrow[r] & R .
        \end{tikzcd}
    \end{equation}
    Using the fact that $\iota$ is a monomorphism (\cite[part (1) of Lemma 5.10]{bongiorno1}), it is easy to see that the datum of $X$, $\hat{R}$ together with $\hat{s}, \hat{t}, \hat{e}, \hat{i}, \hat{c}$ yields a formal groupoid $\hat{\mathcal{X}}$.
    Since $\hat{e}$ is a thickening, $\hat{\mathcal{X}}$ is in fact an infinitesimal groupoid.
    By \cite[part (2) and (3) of Lemma 5.10]{bongiorno1}, $\iota$ is locally of formal finite presentation and formally \'{e}tale.
    Hence $\hat{s}$ is locally of formal finite presentation and, if $s$ is formally smooth, so is  $\hat{s}$.
    We deduce that $\hat{\mathcal{X}} = \hat{R} \rightrightarrows X$ is an infinitesimal groupoid locally of formal finite presentation and it is formally smooth whenever $\mathcal{X}$ is.
    Furthermore, the morphism $\iota$ induces a morphism of formal groupoids $\iota : \hat{\mathcal{X}} \rightarrow \mathcal{X}$.
\end{construction}

The infinitesimal groupoid $\hat{\mathcal{X}}$ has the following universal property: for all infinitesimal groupoids $\mathcal{Z}$ over $S$ and for all morphisms of formal groupoids $g : \mathcal{Z} \rightarrow \mathcal{X}$ over $S$, there exists a unique morphism $\hat{g} : \mathcal{Z} \rightarrow \hat{\mathcal{X}}$ such that $g = \iota \circ \hat{g}$.
This can be checked using the universal property of infinitesimal neighbourhoods.

\begin{example}
\label{ex:de_rham_stack}
    Let $X$ be a smooth scheme over a field $k$.
    The \emph{de Rham stack} $X_{\mathrm{dR}}$ is defined as (the stackification of) the functor
    \begin{align}
    \label{eq:de_rham_stack}
        X_{\mathrm{dR}} : \mathrm{Sch}/k &\longrightarrow \mathrm{Sets} \\
        T &\longrightarrow \mathrm{Hom}_k \left( T_{\mathrm{red}}, X \right). \nonumber
    \end{align}
    This object is used to define crystalline cohomology.
    This should be thought as the infinitesimal neighbourhood of the diagonal.
    To see this, consider the groupoid $X \times_k X$ on $X$.
    Its associated infinitesimal groupoid is $\left( X \times_k X \right)_{\widehat{\Delta}}$, the completion of the diagonal morphism $\Delta_{X/k} : X \rightarrow X \times_k X$.
    It is well known that, for any scheme $T$ over $k$, the natural map
    \begin{align}
    \label{eq:de_rham_stack_proof}
        \left[ \mathrm{Hom}_k \left( T, \left( X \times_k X \right)_{\widehat{\Delta}} \right)
        \rightrightarrows
        \mathrm{Hom}_k \left( T, X \right) \right]
        &\longrightarrow
        \mathrm{Hom}_k \left( T_{\mathrm{red}}, X \right) \\
        \left( T \rightarrow X \right)
        &\longrightarrow
        \left(T_{\mathrm{red}} \rightarrow T \rightarrow X \right) \nonumber
    \end{align}
    is an isomorphism of groupoids up to an \'{e}tale cover of $T$ (see \cite[\href{https://stacks.math.columbia.edu/tag/0CHJ}{Lemma 0CHJ}]{stacks-project}).
\end{example}
\section{Infinitesimal equivalence relations}
\label{sec:infinitesimal_equivalence_relations}

In this section, we show that the infinitesimal groupoid associated to a Deligne--Mumford groupoid is an infinitesimal equivalence relation. \par

\begin{definition}
\label{def:equivalence_relation}
    A \emph{formal equivalence relation} is a formal groupoid $R \rightrightarrows X$ such that $j = t \times s$ is a monomorphism.
    When $R$ is in addition an infinitesimal groupoid, then it is an \emph{infinitesimal equivalence relation}.
\end{definition}

\begin{remark}
\label{rem:equivalence_relation}
    By reducing to the case of sets, it is straightforward to see that a formal groupoid $R \rightrightarrows X$ is a formal equivalence relation if and only if the stabiliser $R \times_{\left( j, \Delta_{X/S} \right)} X$ is trivially equal to $X$.
\end{remark}

\begin{definition}
\label{def:formal_dm_stack}
    Let $R$ be a formal groupoid on a formal scheme $X$ over a formal scheme $S$.
    Then $R \rightrightarrows X$ is a \emph{formal Deligne--Mumford groupoid} if $s : R \rightarrow X$ is formally smooth and locally of formal finite presentation and $j = t \times s : R \rightarrow X \times_S X$ is formally unramified.
    If furthermore $X$ is locally of formal finite presentation over $S$, then $R \rightrightarrows X$ is a formal Deligne--Mumford groupoid \emph{locally of formal finite presentation} over $S$.
\end{definition}

\begin{lemma}
\label{lem:unramified_diagonal_infinitesimal}
    Let $R \rightrightarrows X$ be a formal groupoid over a formal scheme $S$.
    Suppose that $s : R \rightarrow X$ is locally of formal finite presentation and that $j : R \rightarrow X \times_S X$ is formally unramified.
    Then the infinitesimal groupoid $\hat{R}$ associated to $R$ is an infinitesimal equivalence relation.
    In particular the infinitesimal groupoid associated to a Deligne--Mumford groupoid is a formally smooth infinitesimal equivalence relation.
\end{lemma}

\begin{proof}
    Since $s$ is locally of formal finite presentation, we can consider the associated infinitesimal groupoid $\hat{R} \rightrightarrows X$.
    Let $P$ be the stabiliser of $R \rightrightarrows X$ and consider the following commutative diagram
    \begin{equation}
    \label{diag:stabiliser_trivial}
        \begin{tikzcd}
            X \arrow[d, "\mathds{1}"] \arrow[r] & \hat{P} \arrow[d] \arrow[r] & P \arrow[d] \arrow[r] & X \arrow[d, "\Delta_{X/S}"] \\
            X \arrow[r, "\hat{e}"] & \hat{R} \arrow[r] & R \arrow[r, "j"] \arrow[r] & X \times_S X.
        \end{tikzcd}
    \end{equation}
    We note that every square is Cartesian.
    Indeed, the right-most square is Cartesian by definition and the total rectangle is obviously Cartesian.
    Then the pasting Lemma for pull-backs shows that the combination of the left-most and centre squares is Cartesian.
    Finally, \cite[Lemma 5.7]{bongiorno1} shows that both the left-most and centre squares are Cartesian.
    In particular, this shows that the infinitesimal neighbourhood $\hat{P}$ of the stabiliser $P$ of $R \rightrightarrows X$ is the stabiliser of $\hat{R} \rightrightarrows X$. \par

    We now show that the infinitesimal neighbourhood of $X \rightarrow P$ is $X$.
    This immediately shows that $X = \hat{P}$ so that $\hat{R}$ is an infinitesimal equivalence relation.
    By construction, $X \rightarrow P$ is a section of $P \rightarrow X$, which is by assumption formally unramified.
    Suppose that
    \begin{equation}
    \label{diag:infinitesimal_neighbourhood_unramified_i}
        \begin{tikzcd}
            T \arrow[d] \arrow[r] & T^{\prime} \arrow[d] \arrow[ld, dashed]\\
            X \arrow[r] & P.
        \end{tikzcd}
    \end{equation}
    is a commutative solid square where $T \rightarrow T^{\prime}$ is a thickening of affine schemes with square-zero kernel ideal.
    We want to find a unique compatible dashed morphism.
    Composing $T^{\prime} \rightarrow P$ with $P \rightarrow X$ gives a morphism $T^{\prime} \rightarrow X$.
    The upper triangle of Diagram (\ref{diag:infinitesimal_neighbourhood_unramified_i}) is clearly commutative, hence we only have to check that the lower triangle is commutative.
    This amounts to checking that the two morphisms $T^{\prime} \rightrightarrows P$ are equal. 
    By construction, both morphisms fit in the commutative diagram
    \begin{equation}
    \label{diag:infinitesimal_neighbourhood_unramified_ii}
        \begin{tikzcd}
            T \arrow[d] \arrow[r] & T^{\prime} \arrow[d] \arrow[dl, shift left = 0.5ex] \arrow[dl, shift right = 0.5ex] \\
            P \arrow[r] & X.
        \end{tikzcd}
    \end{equation}
    Since $P \rightarrow X$ is formally unramified, there is at most one morphism $T^{\prime} \rightarrow P$.
    This shows that Diagram (\ref{diag:infinitesimal_neighbourhood_unramified_i}) is commutative.
    Furthermore, since $X \rightarrow P$ is a monomorphism, the dashed we have constructed in Diagram (\ref{diag:infinitesimal_neighbourhood_unramified_i}) is unique.
    Hence we have shown that $X$ has the universal property of the infinitesimal neighbourhood of $X \rightarrow P$ for all thickenings of affine schemes with square-zero kernel ideal.
    In fact, it is sufficient to check for this special type of thickenings (\cite[Lemma 5.12]{bongiorno1}), hence $X$ is the infinitesimal neighbourhood of $X \rightarrow P$.
\end{proof}

Next we show that the diagonal $j$ of an infinitesimal equivalence relation is a closed immersion.

\begin{lemma}
\label{lem:monomorphism_closed}
    Let $R \rightrightarrows X$ be an infinitesimal equivalence relation locally of formal finite presentation over a formal scheme $S$ and let $\left( X \times_S X \right)_{\widehat{\Delta}}$ denote the infinitesimal neighbourhood of the diagonal immersion $\Delta_{X/S}$.
    Then there is an induced morphism
    \begin{align}
    \label{eq:monomorphism_closed}
        \hat{j} : R \rightarrow \left( X \times_S X \right)_{\widehat{\Delta}}
    \end{align}
    which is a closed immersion of formal schemes.
\end{lemma}

\begin{proof}
    The assumptions imply that $\Delta_{X/S}$ is locally of formal finite presentation (\cite[Lemma 4.13]{bongiorno1}).
    Then existence of $\hat{j}$ follows by applying the universal property of the infinitesimal neighbourhood of $\Delta_{X/S}$ to the thickening $e : X \rightarrow R$. \par

    Since $R$ is an equivalence relation, the diagram
    \begin{equation}
    \label{diag:monomorphism_closed}
    	\begin{tikzcd}
    		X \arrow[d, "\mathds{1}"] \arrow[r, "e"] & R \arrow[d, "\hat{j}"] \\
    		X \arrow[r, "\hat{\Delta}_{X/S}"] & \left( X \times_S X \right)_{\widehat{\Delta}}.
    	\end{tikzcd}
    \end{equation}
    is Cartesian, where $\hat{\Delta}_{X/S}$ is a thickening of formal schemes locally of formal finite presentation.
    Now the result follows from \cite[Lemma 4.16]{bongiorno1}.
\end{proof}

\begin{comment}
    \begin{lemma}
    \label{lem:regular_immersion_smoothness}
        Let $R \rightrightarrows X$ be an infinitesimal groupoid over a formal scheme $S$ and let $\mathscr{K}$ be the kernel ideal of $\mathscr{O}_R$ corresponding to the thickening $e$.
        Suppose that $s : R \rightarrow X$ is formally smooth, then $\mathscr{K} / \mathscr{K}^2$ is a locally formally projective sheaf on $X$ and the natural morphism of graded topological sheaves on $X$
        \begin{align}
        \label{eq:symmetric_regular_thickening}
            \mathrm{Sym}_X^{\bullet} \, \mathscr{K} / \mathscr{K}^2 \rightarrow \oplus_{n \in \mathbb{N}} \, \mathscr{K}^n / \mathscr{K}^{n+1}
        \end{align}
        is an isomorphism.
    \end{lemma}
    
    \begin{proof}
        \typeout{WARNING: finish proof}
        Use \cite[Theorem 19.5.3, page 91]{MR0173675}.
    \end{proof}
    
    In fact, if being formally smooth were shown to be a local property for formal schemes, then the converse would also hold.
\end{comment}
\section{Infinitesimal restrictions}
\label{sec:restrictions}

In this section, we define the infinitesimal restriction of a groupoids and compare this notion to usual restriction and completion along a closed point.

\begin{definition}
\label{def:restriction}
    Let $\mathcal{X} = R \rightrightarrows X$ be a formal groupoid over a formal scheme $S$ and let $g : Z \rightarrow X$ be a morphism of formal schemes over $S$.
    The \emph{restriction} of $R$ via $g$, denoted by $\mathcal{X}|_Z = R|_Z \rightrightarrows Z$, is defined as the fibre product in the diagram
    \begin{equation}
    \label{diag:restriction_definition}
        \begin{tikzcd}
            R|_Z \arrow[d] \arrow[r] & Z \times_S Z \arrow[d, "g \times g"] \\
            R \arrow[r, "j"] & X \times_S X.
        \end{tikzcd}
    \end{equation}
    This is a formal groupoid on $Z$ over $S$.
    If the source morphism $s : R \rightarrow X$ and $g : Z \rightarrow X$ are morphisms locally of formal finite presentation, the \emph{infinitesimal restriction} of $\mathcal{X}$ via $g$, denoted by $\hat{\mathcal{X}}|_Z = R\hat{|}_Z \rightrightarrows Z$, is the infinitesimal groupoid associated to $R|_Z \rightrightarrows Z$.
\end{definition}

\begin{remark}
\label{rem:restriction_keel_mori}
    There is an alternative way to define the restriction of a formal groupoid: it can be constructed by the following iterated fibre products:
    \begin{equation}
    \label{diag:restriction}
        \begin{tikzcd}
            R|_Z \arrow[r, "g_s^{\prime}"] \arrow[d, "g_t^{\prime}"] & R \times_{(s, g)} Z \arrow[r, "s^{\prime}"] \arrow[d, "g_t"] & Z \arrow[d, "g"] \\
            Z \times_{(g, t)} R \arrow[r, "g_s"] \arrow[d, "t^{\prime}"] & R \arrow[r, "s"] \arrow[d, "t"'] & X \\
            Z \arrow[r, "g"] & X.
        \end{tikzcd}
    \end{equation}
    This fact is observed in \cite[Remark 2.6]{MR1432041}.
\end{remark}

\begin{remark}
\label{rem:infinitesimal_restriction_defined}
    Note that, when defining infinitesimal restriction, we require the assumptions of formal finite presentation in order to establish that $R|_Z$ is a formal groupoid locally of formal finite presentation and apply Construction \ref{cons:infinitesimal_from_formal}.
    Indeed, the source morphism $s_Z : R|_Z \rightarrow Z$ is simply the morphism $s^{\prime} \circ g_s^{\prime} : R|_Z \rightarrow Z$ from Diagram (\ref{diag:restriction}).
    Since both $g$ and $s$ are locally of formal finite presentation, so is $s_Z$. 
    Therefore we can construct the infinitesimal groupoid associated to $R|_Z$.
\end{remark}

\begin{remark}
\label{rem:infinitesimal_restriction_monomorphism}
    If $R \rightrightarrows X$ is an infinitesimal groupoid over $S$ whose source morphism is locally of finite presentation and $g : Z \rightarrow X$ is a monomorphism which is locally of formal finite presentation, then $R \hat{|}_Z = R |_Z$.
    This follows by observing that the commutative diagram
    \begin{equation}
    \label{diag:infinitesimal_restriction_monomorphism}
        \begin{tikzcd}
            Z \arrow[d, "g"] \arrow[r, "e_Z"] & R|_Z \arrow[d] \arrow[r, "j_Z"] & Z \times_S Z \arrow[d, "g \times g"] \\
            X \arrow[r, "e"] & R \arrow[r, "j"] & X \times_S X
        \end{tikzcd}
    \end{equation}
    is Cartesian and that the base change of a thickening is a thickening (\cite[Lemma 3.7]{bongiorno1}).
\end{remark}

\begin{assumptions}
\label{ass:restrictions}
    We consider the following set-up: $R \rightrightarrows X$ is a formal groupoid over a formal scheme $S$ whose source morphism $s : R \rightarrow X$ is locally of formal finite presentation and $x \in X$ is a closed point whose associated closed immersion $\spec \kappa(x) \rightarrow X$ is locally of formal finite presentation.
    We let $\hat{X}$ denote the infinitesimal neighbourhood of $\spec \kappa(x) \rightarrow X$ (see Remark \ref{rem:infinitesimal_restriction_defined}) and let $\hat{R}_x$ denote the infinitesimal neighbourhood of the composition $\spec \kappa(x) \rightarrow X \xrightarrow{e} R$.
\end{assumptions}

In the next lemma we show that infinitesimally restricting a formal groupoid to a closed point $x$ is the same as completing $R$ along $x$.

\begin{lemma}
\label{lem:formal_restriction}
    Under Assumptions \ref{ass:restrictions}, we have a natural isomorphism
    \begin{align}
    \label{eq:formal_restriction}
        R\hat{|}_{\hat{X}} \xrightarrow{\sim} \hat{R}_x.
    \end{align}
\end{lemma}

\begin{proof}
    We first construct a morphism $R\hat{|}_{\hat{X}} \rightarrow \hat{R}_x$.
    Note that $\spec \kappa(x) \rightarrow \hat{X} \rightarrow R\hat{|}_{\hat{X}}$ is a composition of thickenings. Applying the universal property of infinitesimal neighbourhoods to the diagram
    \begin{equation}
    \label{diag:formal_restriction}
        \begin{tikzcd}
            \spec \kappa(x) \arrow[d, "\mathds{1}"] \arrow[rr] & & R\hat{|}_{\hat{X}} \arrow[d] \arrow[dl, dashed] \\
            \spec \kappa(x) \arrow[r] & \hat{R}_x \arrow[r] & R
        \end{tikzcd}
    \end{equation}
    yields a unique compatible morphism. \par

    Now we construct the inverse morphism.
    Note that the infinitesimal neighbourhood of the closed immersion $\spec \kappa(x) \rightarrow X \times_S X$ is $\hat{X} \times_S \hat{X}$.
    Indeed, the universal property can be verified by applying \cite[Lemma 5.7]{bongiorno1} to the base change of the infinitesimal neighbourhood
    \begin{align}
    \label{eq:infinitesimal_neighbourhood_product}
        \spec \kappa(x) \rightarrow \hat{X} \rightarrow X
    \end{align}
    by the morphism $X \times_S X \rightarrow X$ and using the fact that the commutative diagram
    \begin{equation}
    \label{diag:infinitesimal_products}
        \begin{tikzcd}
            \hat{X} \times_S \hat{X} \arrow[d] \arrow[r] & X \times_S \hat{X} \arrow[d] \\
            \hat{X} \times_S X \arrow[r] & X \times_S X
        \end{tikzcd}
    \end{equation}
    is Cartesian.
    Therefore the morphism $\hat{R}_x \rightarrow X \times_S X$ induces a unique morphism $\hat{R}_x \rightarrow \hat{X} \times_S \hat{X}$.
    By the fibre product in Diagram (\ref{diag:restriction_definition}), there exists a unique morphism $\hat{R}_x \rightarrow R|_{\hat{X}}$.
    Finally, since $\spec \kappa(x) \rightarrow \hat{R}_x$ is a thickening, there exists a unique morphism $\hat{R}_x \rightarrow R\hat{|}_{\hat{X}}$. \par

    The morphisms we have constructed are unique over $R$.
    Since both $\hat{R}_x \rightarrow R$ and $R\hat{|}_{\hat{X}} \rightarrow R$ are monomorphisms (\cite[part (1) of Lemma 5.10]{bongiorno1}), we conclude that (\ref{eq:formal_restriction}) holds.
\end{proof}

In the next lemma we show that infinitesimally restricting a formal groupoid is the same as infinitesimally restricting its associated infintiesimal groupoid.

\begin{lemma}
\label{lem:formal_formal_restriction}
    Under Assumptions \ref{ass:restrictions}, let $\hat{R}$ be the infinitesimal groupoid associated to $R$ from Construction \ref{cons:infinitesimal_from_formal}.
    Then there exists a natural isomorphism
    \begin{align}
    \label{eq:formal_formal_restriction}
        \hat{R} \hat{|}_{\hat{X}} \xrightarrow{\sim} R \hat{|}_{\hat{X}}
    \end{align}
\end{lemma}

\begin{proof}
    By Lemma \ref{lem:formal_restriction}, it suffices to show that the infinitesimal neighbourhoods of $x$ in $\hat{R}$ and $R$ are naturally isomorphic.
    This follows from \cite[Lemma 5.13]{bongiorno1} once we observe that $\hat{R} \rightarrow R$ is formally \'{e}tale (\cite[part (3) of Lemma 5.10]{bongiorno1}).
\end{proof}

\begin{comment}
    Infinitesimal restriction can be characterised by the following universal property: for all infinitesimal groupoids $\mathcal{Z}$ over $S$ and for all morphisms of groupoids $Z \rightarrow \mathcal{Z}$ and $\mathcal{Z} \rightarrow \mathcal{X}$ over $S$ such that the solid diagram
    \begin{equation}
    \label{diag:restriction_universal_property}
        \begin{tikzcd}
            Z \arrow[d] \arrow[dr] \arrow[rr, "g"] & & X \arrow[d] \\
            \mathcal{Z} \arrow[r, dashed] \arrow[rr, bend right] & \hat{\mathcal{X}}|_Z \arrow[r] & \mathcal{X},
        \end{tikzcd}
    \end{equation}
    is commutative, there exists a unique compatible dashed morphism.
    This can be checked combining the universal property of restriction of groupoids and the universal property of the infinitesimal groupoids associated to formal groupoids.
\end{comment}
\section{Transversals}
\label{sec:transversals}

In this section, we define what is a transversal (or slice) of a formal groupoid and we show it always exists under Noetherian assumptions.
Our notion of transversal is fairly strict and it only makes sense for formal Deligne--Mumford groupoids.
As a result, we restrict to this particular case.

\begin{definition}
\label{def:transversal}
    Let $R \rightrightarrows X$ be a formal Deligne--Mumford groupoid over a formal scheme $S$ and let $x \in X$ be a closed point whose associated closed immersion $\spec \kappa(x) \rightarrow X$ is locally of formal finite presentation.
    A \emph{transversal} of $R \rightrightarrows X$ through $x$ is an immersion of formal schemes $g : W \rightarrow X$ locally of formal finite presentation such that $x \in \mathrm{im}\,(g)$ and the induced morphism
    \begin{align}
    \label{eq:transveral_condition}
        p : N := W \times_{(g, t)} R \xrightarrow{g_s} R \xrightarrow{s} X
    \end{align}
    from Diagram (\ref{diag:restriction}) induces a formal-local isomorphism at $x$.
    More precisely, let $\hat{X}$ and $\hat{N}$ be the infinitesimal neighbourhoods of $x \in X$ and $(e(x), x) \in N$ respectively, then $p : N \rightarrow X$ induces an isomorphism $\hat{p} : \hat{N} \xrightarrow{\sim} \hat{X}$.
\end{definition}

\begin{remark}
\label{rem:neighbourhood_torsor}
    We briefly note that the morphism
    \begin{align}
    \label{eq:neighbourhood_torsor}
        \spec \kappa((e(x),x)) \rightarrow N
    \end{align}
    is an immersion locally of formal finite presentation.
    This ensures existence of $\hat{N}$.
    To this end, we observe that the composition $\spec \kappa((e(x), x) \rightarrow N \xrightarrow{p} X$ is, by assumption, an immersion locally of formal finite presentation.
    Hence, by \cite[Lemma 2.9]{bongiorno1} and \cite[Lemma 4.13]{bongiorno1}, it suffices to show that $p$ is locally of formal finite presentation.
    This is true since $p$ is the composition of two morphisms locally of formal finite presentation.
\end{remark}

The next lemma shows that constructing a transversal only depends on the formal-local structure of a groupoid.

\begin{lemma}
\label{lem:transversal_formal_local}
    Notation as in Definition \ref{def:transversal}, let $\hat{R}$ be the infinitesimal groupoid associated to $R$.
    Then $g : W \rightarrow X$ is a transversal of $R$ if and only if it is a transversal of $\hat{R}$.
\end{lemma}

\begin{proof}
    We simply have to check that the infinitesimal neighbourhoods of $(e(x), x)$ in $W \times_{(g, \hat{t})} \hat{R}$ and $W \times_{(g, t)} R$ are naturally isomorphic.
    By \cite[part (3) of Lemma 5.10]{bongiorno1}, $\hat{R} \rightarrow R$ is formally \'{e}tale.
    Hence the base change $W \times_{(g, \hat{t})} \hat{R} \rightarrow W \times_{(g, t)} R$ is formally \'{e}tale.
    Now the lemma follows from \cite[Lemma 5.13]{bongiorno1}.
\end{proof}

\begin{assumptions}
\label{ass:transversals}
    We consider the following set-up: $R \rightrightarrows X$ is a formal Deligne--Mumford groupoid locally of formal finite presentation over a locally Noetherian formal scheme $S$ and $x \in X$ is a closed point.
    In this case, since $X$ is locally Noetherian (\cite[Lemma 4.14]{bongiorno1}), the closed immersion $l : \spec \kappa(x) \rightarrow X$ is locally of formal finite presentation.
\end{assumptions}

The following construction is based on \cite[Step 2 of Lemma 3.3]{MR1432041}.
We will simply construct a candidate locally closed formal subscheme $W \subseteq X$ and we will only show later in Proposition \ref{prop:existence_transversal} that it is a transversal.

\begin{construction}
\label{cons:transversal}
    Under Assumptions \ref{ass:transversals}, we will construct a transversal of $R \rightrightarrows X$ through $x \in X$. \par
    
    Up to replacing $X$ with an affine neighbourhood of $x$, we may assume $X = \spf A$ is affine.
    Let $\mathfrak{m}$ be the maximal ideal of $A$ corresponding to the closed point $x$.
    Let
    \begin{align}
    \label{eq:completion_point}
        l : \spec \kappa(x) \xrightarrow{\hat{l}} \hat{X} \rightarrow X
    \end{align}
    be the infinitesimal neighbourhood.
    By \cite[Proposition 5.2]{bongiorno1}, $\hat{X} = \spf \hat{A}$, where $\hat{A}$ is the $\mathfrak{m}$-adic completion of $A$.
    Let $\hat{\mathfrak{m}} = \mathfrak{m} \cdot A$ be the maximal ideal of $\hat{A}$ cutting out $x$. \par

    Let $\hat{R}$ be the infinitesimal groupoid associated to $R$ and let 
    \begin{align}
    \label{eq:definition_leaf}
        l^{\prime} : L := \hat{R} \times_{(\hat{s}, \hat{l})} \spec \kappa(x) \rightarrow \hat{R}
    \end{align}
    be the \emph{leaf} through $x$.
    We show that $r := \hat{t} \circ l^{\prime} : L \rightarrow \hat{X}$ is a closed immersion of formal schemes.
    This is the composition of the two vertical arrows in the second column of Diagram (\ref{diag:restriction}) where $W$ and $R$ are replaced by $\spec \kappa(x)$ and $\hat{R}$ respectively. \par

    We first show that $r$ is a closed immersion of affine Noetherian local formal schemes.
    Consider the commutative diagram
    \begin{equation}
    \label{diag:leaf_closed}
        \begin{tikzcd}
            \spec \kappa (x) \arrow[d, "\mathds{1}"] \arrow[r] & L \arrow[d, "r"] \\
            \spec \kappa (x) \arrow[r, "\hat{l}"] & \hat{X}.
        \end{tikzcd}
    \end{equation}
    By Lemma \ref{lem:unramified_diagonal_infinitesimal}, $\hat{R}$ is an infinitesimal equivalence relation, hence Diagram (\ref{diag:leaf_closed}) is Cartesian.
    Since $\hat{l}$ is a thickening of formal schemes locally of formal finite presentation, \cite[Lemma 4.16]{bongiorno1} implies that $r$ is a closed immersion of formal schemes.
    Since $\hat{X}$ is affine, so is $L$.
    Let $L = \spf B$ and let $\mathfrak{n} = \mathfrak{m} \cdot B$ be its maximal ideal cutting out $x$.
    Since $r$ is adic, $\mathfrak{n}$ is an ideal of definition of $B$.
    Now $A$ is Noetherian by assumption, hence $\hat{A}$ is an adic Noetherian local ring.
    Since $r$ is a closed immersion, $B$ is also an adic Noetherian local ring and $r$ induces a surjective morphism of adic Noetherian local rings. \par

    Now we construct $W \subseteq X$.
    Let $(\bar{f}_1, \ldots, \bar{f}_d)$ be a minimal set of generators of the maximal ideal $\mathfrak{n} \subseteq B$.
    We have that $\mathfrak{n} = \hat{\mathfrak{m}} \cdot B = \mathfrak{m} \cdot B$, hence the elements $\{ f_i \}_{i \leq d}$ can be assumed to be in the image (i.e. not in the larger ideal extension) of $\mathfrak{m}$ in $B$.
    Now lift such elements to elements $(f_1, \ldots, f_d) \subseteq \mathfrak{m}$ and define $W \subseteq \spf A$ to be the closed subscheme cut out by $(f_1, \ldots, f_d)$.
    By construction, $W \subseteq X$ is an immersion of formal schemes locally of formal finite presentation such that
    \begin{align}
    \label{eq:orthogonal_intersection}
        \hat{W} \times_{\hat{X}} L = \spec \kappa (x),
    \end{align}
    where $\hat{W}$ is the infinitesimal neighbourhood of $x \in W$.
\end{construction}

In fact, since $R$ is formally smooth, the minimal set of generators in the previous construction must be a regular sequence.

\begin{lemma}
\label{lem:leaf_smooth}
    Under Assumptions \ref{ass:transversals} and using the notation of Construction \ref{cons:transversal}, $(\bar{f}_1, \ldots, \bar{f}_d)$ is a regular system of parameters of $\mathfrak{n} \subseteq B$.
\end{lemma}

\begin{proof}
    Since $s : \hat{R} \rightarrow X$ is formally smooth, so is its base change
    \begin{align}
    \label{eq:leaf_smooth}
        L = \hat{R} \times_{(\hat{s}, \hat{l})} \spec \kappa(x) \rightarrow \spec \kappa(x).
    \end{align}
    Recall from Construction \ref{cons:transversal} that $L = \spf B$ is an adic Noetherian local ring with maximal ideal $\mathfrak{n}$ and residue field $\kappa(x)$.
    Now \cite[\href{https://stacks.math.columbia.edu/tag/07EI}{Lemma 07EI}]{stacks-project} implies that $B$ is a regular local ring and \cite[\href{https://stacks.math.columbia.edu/tag/00NQ}{Lemma 00NQ}]{stacks-project} implies that $(\bar{f}_1, \ldots, \bar{f}_d)$ is a \emph{regular} system of parameters.
\end{proof}

We finally show that the locally closed formal subscheme $W \subseteq X$ is a transversal.

\begin{proposition}
\label{prop:existence_transversal}
    Under Assumptions \ref{ass:transversals}, the immersion of formal schemes $g : W \rightarrow X$ of Construction \ref{cons:transversal} is a transversal through $x$, i.e. the morphism
    \begin{align}
        \hat{p} : \hat{N} \xrightarrow{\,\sim\,} \hat{X}
    \end{align}
    constructed in Definition \ref{def:transversal} is an isomorphism.
\end{proposition}

The proof of Proposition \ref{prop:existence_transversal} is articulated into six steps:
\begin{enumerate}
    \item We reduce to the case where $R$ is an infinitesimal equivalence relation.
    \item We show that $\hat{p}$ induces a morphism of adic Noetherian local rings.
    \item We show that $\hat{p}$ is flat as a morphism of discrete rings.
    \item We show that $\hat{p}$ is a formally smooth morphism.
    \item We show that $\hat{p}$ is a thickening of formal schemes.
    \item We show that $\hat{p}$ is an isomorphism.
\end{enumerate}
%The proof of Step 3 is based on \cite[Step 2 of Lemma 3.3]{MR1432041}.

\begin{proof}
[Proof of Proposition \ref{prop:existence_transversal}]
    \emph{Step 1.}
    By Lemma \ref{lem:transversal_formal_local}, when verifying whether $g : W \rightarrow X$ is a transversal, $R$ may be replaced with its associated infinitesimal groupoid $\hat{R}$.
    Since $R$ is a formal Deligne--Mumford groupoid, $\hat{R}$ is an infinitesimal equivalence relation whose source morphism is formally smooth and locally of formal finite presentation (Lemma \ref{lem:unramified_diagonal_infinitesimal}).
    As a result, there is no harm in further assuming $R$ to be an equivalence relation. \par

    \emph{Step 2.}
    By assumption $R$ is locally of formal finite presentation over a locally Noetherian scheme $S$, hence $R$ is locally Noetherian (\cite[Lemma 4.14]{bongiorno1}).
    Furthermore, we know that $X = \spf A$ and $\hat{X} = \spf \hat{A}$.
    It follows that $R = \spf \Gamma$ is affine (\cite[Lemma 3.8]{bongiorno1}) and Noetherian.
    Let $\hat{R}_x$ be the infinitesimal neighbourhood of $\spec \kappa(e(x)) \in R$.
    We see that $\hat{R}_x = \spf \hat{\Gamma}$, where $\hat{\Gamma}$ is the adic completion of $\Gamma$ at the maximal ideal corresponding to $e(x)$.
    Therefore $\hat{\Gamma}$ is an adic Noetherian local ring.
    It follows that $\hat{s} : \hat{R}_x \rightarrow \hat{X}$ induces a morphism $\hat{s}^{\#} : \hat{A} \rightarrow \hat{\Gamma}$ of adic Noetherian local rings.
    Therefore, in order to show that $\hat{p}$ is a morphism of adic Noetherian local rings, it suffices to show that so is $\hat{g}_s : \hat{N} \rightarrow \hat{R}_x$.
    To this end, let $t^{\#} : A \rightarrow \Gamma$ denote the morphism of adic rings associated to $t : R \rightarrow X$ and let $h_i = t^{\#}(f_i)$ for all $i \leq d$.
    By construction, $W = \spf \left( A / (f_1, \ldots, f_d) \right)$.
    Since $g$ is a closed immersion, so is $g_s$, hence $N = \spf \left( \Gamma / (h_1, \ldots, h_d) \right)$.
    Let $\hat{W}$ be the infinitesimal neighbourhood of $\spec \kappa(x) \rightarrow W$.
    Combining \cite[Lemma 5.7]{bongiorno1} together with the fact that $g$ is a monomorphism implies that $\hat{W} = W \times_X \hat{X}$.
    Similarly, since $g_s$ is a monomorphism, $\hat{N} = N \times_R \hat{R}_x$. This yields the following affine descriptions
    \begin{align}
    \label{eq:infinitesimal_torsor_i}
        \hat{W} &= \spf \left( \hat{A} / (f_1, \ldots, f_d) \right), \\
    \label{eq:infinitesimal_torsor_ii}
        \hat{N} &= \spf \left( \hat{\Gamma} / (h_1, \ldots, h_d) \right).
    \end{align}
    It follows that $\hat{p} = \hat{s} \circ \hat{g}_s : \hat{N} \rightarrow \hat{X}$ is a morphism of adic Noetherian local rings. \par
    
    \emph{Step 3.}
    We note that the morphism $l^{\prime} : L \rightarrow R$ in (\ref{eq:definition_leaf}) factors through $\hat{l}^{\prime} : L \rightarrow \hat{R}_x$.
    This follows from the universal property of infinitesimal neighbourhoods applied to the thickening $\spec \kappa(x) \rightarrow L$.
    Furthermore, since $\hat{R}_x \rightarrow R$ is a monomorphism (\cite[part (1) of Lemma 5.10]{bongiorno1}), $L = \hat{R}_x \times_{(\hat{s}, \hat{l})} \spec \kappa(x)$.
    Together with (\ref{eq:orthogonal_intersection}), this yields the following commutative diagram
    \begin{equation}
    \label{diag:flatness_helper}
        \begin{tikzcd}
            \spec \kappa(x) \arrow[d] \arrow[r] & L \arrow[d, "\hat{l}^{\prime}"] \arrow[r] & \spec \kappa(x) \arrow[d, "\hat{l}"] \\
            \hat{N} \arrow[d] \arrow[r, "\hat{g}_s"] & \hat{R}_x \arrow[d, "\hat{t}"] \arrow[r, "\hat{s}"] & \hat{X} \\
            \hat{W} \arrow[r, "\hat{g}"] & \hat{X}, &
        \end{tikzcd}
    \end{equation}
    where every square is Cartesian.
    Using (\ref{eq:infinitesimal_torsor_i}) and (\ref{eq:infinitesimal_torsor_ii}), we deduce that Diagram (\ref{diag:flatness_helper}) has the following ring-theoretic counterpart.
    \begin{equation}
    \label{diag:flatness_helper_affine}
        \begin{tikzcd}
            \dfrac{\hat{A}}{\hat{\mathfrak{m}}} & B \arrow[l] & \dfrac{\hat{A}}{\hat{\mathfrak{m}}} \arrow[l] \\
            \dfrac{\hat{\Gamma}}{\left(h_1, \ldots, h_d \right)} \arrow[u] & \hat{\Gamma} \arrow[l] \arrow[u] & \hat{A} \arrow[l, "\hat{s}^{\#}"] \arrow[u] \\
            \dfrac{\hat{A}}{ \left(f_1, \ldots, f_d \right)} \arrow[u] & \hat{A} \arrow[l] \arrow[u, "\hat{t}^{\#}"]. &
        \end{tikzcd}
    \end{equation}
    We denote by $\psi : A \rightarrow \hat{\Gamma} \rightarrow \hat{\Gamma} / (h_1, \ldots, h_d)$ the composition of the two horizontal arrows in the second row.
    This morphism is the ring-theoretic counterpart of $\hat{p} : \hat{N} \rightarrow \hat{X}$. \par

    Let $\psi_k : A \rightarrow \hat{\Gamma} / (h_1, \ldots, h_k)$ for $k \leq d$, where $\psi_0$ is understood to be $\hat{s}^{\#}$.
    Since $\psi_0$ is a formally smooth morphism of local rings, where $\hat{A}$ and $\hat{\Gamma}$ are endowed with the maximal ideal adic topology, \cite[Th\'{e}or\`{e}me 19.7.1, page 200]{MR0173675} implies that $\psi_0$ is flat as a morphism of discrete Noetherian local rings.
    Now consider
    \begin{align}
    \label{eq:multiplication_map}
        \cdot h_1 : \hat{\Gamma} \rightarrow \hat{\Gamma}
    \end{align}
    as a morphism of $\hat{A}$-modules.
    We use \cite[(10.2.4), page 363]{MR0217085} to show that $\psi_1$ is flat.
    To this end, since $\psi_0$ is flat, it suffices to show that the base change of (\ref{eq:multiplication_map}) by the morphism $\hat{A} \rightarrow \hat{A} / \hat{\mathfrak{m}}$ remains injective.
    By construction, the base change morphism is
    \begin{align}
    \label{eq:multiplication_map_leaf}
        \cdot \bar{f}_1 : B \rightarrow B.
    \end{align}
    This is injective.
    Indeed $\bar{f}_1$ is part of a regular system of parameters of $B$ (Lemma \ref{lem:leaf_smooth}).
    We deduce that $\psi_1$ is injective.
    Now applying induction on $\psi_k$ together with the fact that $(\bar{f}_1, \ldots, \bar{f}_d)$ is a regular system of parameters implies that $\psi = \psi_d$ is flat.
    This shows flatness of $\hat{p}$. \par

    \emph{Step 4.}
    From Step 2 and Step 3, we know that $\hat{p}$ is a morphism of adic Noetherian local rings which is flat as a morphism of discrete rings.
    Therefore, by \cite[Th\'{e}or\`{e}me 19.7.1, page 200]{MR0173675}, showing the $\hat{p}$ is formally smooth amounts to showing that the morphism $\hat{N} \times_{\hat{X}} \spec \kappa(x) \rightarrow \spec \kappa(x)$ is formally smooth.
    In fact, Diagram (\ref{diag:flatness_helper}) immediately shows it is an isomorphism, hence $\hat{p}$ is formally smooth.
    
    \emph{Step 5.}
    Note that the diagram
    \begin{equation}
    \label{diag:transversal_closed_immersion}
        \begin{tikzcd}
            \spec \kappa(x) \arrow[d] \arrow[r, "\mathds{1}"] & \spec \kappa(x) \arrow[d, "\hat{l}"] \\
            \hat{N} \arrow[r, "\hat{p}"] & \hat{X}
        \end{tikzcd}
    \end{equation}
    is Cartesian.
    This was observed in Diagram (\ref{diag:flatness_helper}).
    Since $\hat{l}$ is a thickening of formal schemes locally of formal finite presentation, \cite[Lemma 4.16]{bongiorno1} implies that $\hat{p}$ is a thickening of formal schemes. \par
    
    \emph{Step 6.}
    Clearly, $\hat{p}$ is a morphism locally of finite presentation since so are $\hat{g}$ and $\hat{s}$.
    From Step 4 and Step 5, we know that $\psi$ is a formally smooth thickening of formal schemes, hence it is an isomorphism (\cite[Lemma 4.15]{bongiorno1}).
\end{proof}

\begin{remark}
\label{rem:involution_transversal}
    If $W$ is a transversal of $R \rightrightarrows X$ through $x \in X$, by definition, $s \circ g_s : W \times_{(g,t)} R \rightarrow X$ induces a formal-local isomorphism.
    Composing with the structural involution $i : R \rightarrow R$ yields that
    \begin{align}
    \label{eq:involution_transversal}
        t \circ g_t : R \times_{(s,g)} W \rightarrow X
    \end{align}
    also induces a formal-local isomorphism.
\end{remark}
\section{Frobenius theorem}
\label{sec:frobenius_theorem}

In this section, we prove the main theorems.
The most important work has already been done in the previous section and now it suffices to run some category-theoretic arguments.
\begin{assumptions}
\label{ass:single_point}
    We consider the following set-up: $R \rightrightarrows X$ is an infinitesimal Deligne--Mumford groupoid locally of formal finite presentation over a locally Noetherian formal scheme $S$ and $|X|$, the topological space of $X$, consists of a single closed point $x$.
    In this case, Proposition \ref{prop:existence_transversal} immediately applies to show existence of a transversal $g : W \rightarrow X$ through $x \in X$, which we fix once and for all.
\end{assumptions}

We begin with an easy consequence of existence of transversals.

\begin{lemma}
\label{lem:q_properties}
    Under Assumptions \ref{ass:single_point}, every square of the following commutative diagram 
    \begin{equation}
    \label{diag:q_construction}
        \begin{tikzcd}
            W \arrow[d, "g"] \arrow[r, "g"] & X \arrow[d, "g_t"] \arrow[r, "q"] & W \arrow[d, "g"] \\
            X \arrow[d, "q"] \arrow[r, "g_s"] & R \arrow[d, "t"] \arrow[r, "s"] & X \\
            W  \arrow[r, "g"] & X. &
        \end{tikzcd}
    \end{equation}
    is Cartesian, where $s \circ g_s = \mathds{1}$ and $t \circ g_t = \mathds{1}$.
    Furthermore, $q$ is a formally smooth morphism locally of formal finite presentation and is split by $g$, i.e. $q \circ g = \mathds{1}_W$.
\end{lemma}

\begin{proof}
    Existence of a transversal follows immediately from Proposition \ref{prop:existence_transversal}. \par
    
    Let $N = W \times_{(g, t)} R$.
    We show that $x \in N$ is a thickening, so that $N = \hat{N} = \hat{X}$.
    Since $x \in X$ is a thickening and $X \rightarrow R$ is a thickening, we know that $x \in R$ is a thickening.
    But then $N$ is a non-empty formal subscheme of $R$ and thus its topological space is a single point $x$.
    Now \cite[Lemma 2.9]{bongiorno1} implies that $x \in N$ is a thickening.
    This shows that the bottom-left square of Diagram (\ref{diag:q_construction}) is Cartesian.
    Since $W$ is a transversal, $s \circ g_s = \mathds{1}$.
    A specular argument using Remark \ref{rem:involution_transversal} shows that the top-right square is Cartesian and that $t \circ g_t = \mathds{1}$.
    It now follows easily that the remaining top-right square is Cartesian too. \par

    It is straightforward to observe that $q$ is formally smooth and locally of finite presentation as it is the base change of $t$.
    We finally show that $q \circ g = \mathds{1}_W$.
    Note that $R\hat{|}_{W} = W \rightrightarrows W$ (Lemma \ref{lem:formal_restriction}).
    This shows that the composition $q \circ g$ is the target morphism of the groupoid $W \rightrightarrows W$, hence it has to be the identity.
\end{proof}

The next proposition shows that the morphism $q : X \rightarrow W$ is the \emph{effective} geometric quotient of $X$ by $R$.

\begin{proposition}
\label{prop:q_invariance}
    Under Assumptions \ref{ass:single_point}, let $q : X \rightarrow W$ be the $S$-morphism constructed in Lemma \ref{lem:q_properties} and consider the fibre product $X \times_W X$ as a groupoid on $X$.
    Then there exists an isomorphism of groupoids
    \begin{align}
    \label{eq:crucial_isomorphism}
        \Phi : X \times_W X \xrightarrow{\;\sim\;} R
    \end{align}
    on $X$.
\end{proposition}

\begin{proof}
    We first construct $\Phi$.
    Consider the diagram
    \begin{equation}
    \label{diag:composition_phi}
        \begin{tikzcd}
            X \times_W X \arrow[d, "g_t \circ \mathrm{pr}_1"] \arrow[r, "g_s \circ \mathrm{pr}_2"] & R \arrow[d, "t"] \\
            R \arrow[r, "s"] & X.
        \end{tikzcd}
    \end{equation}
    This is commutative.
    Indeed, by construction and Diagram (\ref{diag:q_construction}), both composition of morphisms are given by $X \times_W X \rightarrow W \rightarrow X$.
    This gives a morphism
    \begin{align}
        \Phi : X \times_W X \xrightarrow{g_{s,t}} R \times_{(s, t)} R \xrightarrow{c} R.
    \end{align}
    Furthermore $\Phi$ is a morphism over $X \times_S X$.
    More precisely,
    \begin{align}
    \label{eq:source_compatibility}
        s \circ \Phi &= \mathrm{pr}_2 : X \times_W X \rightarrow X, \\
    \label{eq:target_compatibility}
        t \circ \Phi &= \mathrm{pr}_1 : X \times_W X \rightarrow X.
    \end{align}
    Indeed, using that $t \circ c = t \circ \mathrm{pr}_1$ (\cite[\href{https://stacks.math.columbia.edu/tag/02YE}{Lemma 02YE}]{stacks-project}) and $t \circ g_t = \mathds{1}_X$ yields
    \begin{align}
    \label{eq:linear_over_product}
        t \circ \Phi = t \circ c \circ g_{s,t} = t \circ \mathrm{pr}_1 \circ g_{s,t} = t \circ g_t \circ \mathrm{pr}_1 = \mathrm{pr}_1.
    \end{align}
    A specular argument yields that $s \circ \Phi = \mathrm{pr}_2$. \par
    
    Next we show that $q$ is $R$-invariant, i.e. $q \circ s = q \circ t$.
    Since the bottom-left square of Diagram (\ref{diag:q_construction}) is commutative, we have that $g \circ q = t \circ g_s$.
    Hence $g \circ q \circ s = t \circ g_s \circ s$.
    Since $g_s \circ s = \mathds{1}$, we have that $g \circ q \circ s = t$.
    Now post-composing with $q$ and using $q \circ g = \mathds{1}$ gives that $q \circ s = q \circ t$. \par

    Therefore there is a morphism $\Psi : R \rightarrow X \times_W X$ over $X \times_S X$.
    Since both $R$ and $X \times_W X$ are equivalence relations (Lemma \ref{lem:unramified_diagonal_infinitesimal}) and $\Phi$ and $\Psi$ are morphisms over $X \times_S X$, $\Phi$ and $\Psi$ must be inverse to each other.
    Finally, since any monomorphism into $X \times_S X$ has at most one groupoid structure on $X$, $\Phi$ and $\Psi$ must be morphisms of groupoids.
\end{proof}

The next lemma shows that $q$ is the categorical quotient of $X$ by $R$.
In the case of schemes, this would easily follow from the theory of descent, however this is not yet available for formal schemes.
Instead, we use the fact that $q$ is split by $g$.

\begin{lemma}
\label{lem:q_categorical}
    Under Assumptions \ref{ass:single_point}, the $S$-morphism $q : X \rightarrow W$ constructed in Lemma \ref{lem:q_properties} is the categorical quotient of $X$ by $R$ in the category of formal schemes over $S$.
\end{lemma}

\begin{proof}
    We first recall that $q$ is $R$-invariant as observed in Proposition \ref{prop:q_invariance}. \par
    
    Now, we have to show that for any $R$-invariant $S$-morphism $q^{\prime} : X \rightarrow W^{\prime}$, there exists a unique morphism $w : W \rightarrow W^{\prime}$ such that $q^{\prime} = w \circ q$.
    It is straightforward to see that any compatible morphism $w$ must satisfy $w = q^{\prime} \circ g$, hence if it exists it is unique. \par
    
    We verify that $w = q^{\prime} \circ g$ satisfies $q^{\prime} = w \circ q$.
    Since $q^{\prime}$ is $R$-invariant, the compositions of morphisms in
    \begin{align}
    \label{eq:q_prime_invariant}
        R \xrightrightarrows[t]{s} X \xrightarrow{q^{\prime}} W^{\prime}
    \end{align}
    are equal.
    Define a morphism
    \begin{align}
    \label{eq:twisted_immersion}
        X \xrightarrow{\mathds{1} \times (g \circ q)} X \times_W X \xrightarrow{\Phi} R,
    \end{align}
    where $\Phi$ is the morphism of groupoids constructed in Proposition \ref{prop:q_invariance}.
    Since the compositions in (\ref{eq:q_prime_invariant}) are equal, composing (\ref{eq:twisted_immersion}) with (\ref{eq:q_prime_invariant}) yields two equal morphisms
    \begin{align}
    \label{eq:unique_quotient}
        X \xrightrightarrows[q^{\prime} \circ g \circ q]{q^{\prime}} W^{\prime}.
    \end{align}
    This shows that $q^{\prime} = w \circ q$.
\end{proof}

The next lemma essentially shows that the formal scheme $W$ represents the formal algebraic stack $[X/R]$.
In the statement we try to avoid the word \emph{stack}.

\begin{lemma}
\label{lem:q_representable}
    Under Assumptions \ref{ass:single_point}, for any formal scheme $T$ over $S$, the $S$-morphism $q : X \rightarrow W$ constructed in Lemma \ref{lem:q_properties} induces an isomorphism of groupoids of sets
    \begin{align}
    \label{eq:representable_stack}
        q(T) : \left( \mathrm{Hom}_S \left(T, R \right) \xrightrightarrows[t(T)]{s(T)} \mathrm{Hom}_S \left(T, X \right) \right) \xrightarrow{\;\sim\;} \mathrm{Hom}_S \left(T, W \right).
    \end{align}
\end{lemma}

\begin{proof}
    By Lemma \ref{lem:q_categorical}, $q$ is $R$-invariant, hence $q(T)$ is a morphism of groupoids. \par
    
    We claim that $g(T)$ is the inverse morphism.
    Certainly, since $q \circ g = \mathds{1}_W$, it follows immediately that $q(T) \circ g(T)$ is the identity. \par

    We show that $g(T) \circ q(T)$ is the identity.
    Let $h : T \rightarrow X$ be an $S$-morphism.
    We want to show that $h$ and $g \circ q \circ h$ are equivalent.
    This amounts to finding a morphism $T \rightarrow R$ which gives $h$ when composed with $t$ and $g \circ q \circ h$ when composed with $s$.
    This is readily constructed as
    \begin{align}
    \label{eq:equivalence_morphisms}
        T \xrightarrow{h \times (g \circ q \circ h)} X \times_W X \xrightarrow{\Phi} R,
    \end{align}
    where $\Phi$ is the morphism of groupoids constructed in Proposition \ref{prop:q_invariance}.
\end{proof}

\begin{proof}
[Proof of Theorem \ref{thm:local_structure}]
    By Proposition \ref{prop:q_invariance}, there exists an $S$-morphism $q : X \rightarrow W$ such that $X \times_W X = R$ as groupoids on $X$.
    By Lemma \ref{lem:q_categorical}, $q$ is the categorical quotient in the category of formal schemes over $S$, therefore $q$ and $W$ are unique.
    By Lemma \ref{lem:q_representable}, $q$ represents the stack $[X/R]$.
    By Lemma \ref{lem:q_properties}, $q$ is formally smooth, split and locally of formal finite presentation.
\end{proof}

\begin{proof}
[Proof of Theorem \ref{thm:frobenius}]
    Let $\hat{X}$ be the infinitesimal neighbourhood of the closed point $x \in X$ and let $\hat{R}$ be the infinitesimal groupoid associated to $R$.
    By Lemma \ref{lem:formal_formal_restriction}, $R \hat{|}_{\hat{X}} = \hat{R} \hat{|}_{\hat{X}}$.
    But $\hat{R} \hat{|}_{\hat{X}} \rightrightarrows \hat{X}$ is an infinitesimal Deligne--Mumford groupoid locally of formal finite presentation over a locally Noetherian formal scheme $S$ and the topological space of $\hat{X}$ is a single closed point.
    Therefore Theorem \ref{thm:local_structure} readily applies to give the result.
\end{proof}

\bibliographystyle{alpha}
\bibliography{main}

\end{document}